# FURTHER REFINEMENTS OF GENERALIZED NUMERICAL RADIUS INEQUALITIES FOR HILBERT SPACE OPERATORS

MONIRE HAJMOHAMADI[1], RAHMATOLLAH LASHKARIPOUR[2], MOJTABA BAKHERAD[3]

ABSTRACT. In this paper, we show some refinements of generalized numerical radius inequalities involving the Young and Heinz inequalities. In particular, we present

$$w_p^p(A_1^*T_1B_1, ..., A_n^*T_nB_n)$$
$$\leq \frac{n^{1-\frac{1}{r}}}{2^{\frac{1}{r}}} \Big\| \sum_{i=1}^n [B_i^* f^2(|T_i|)B_i]^{rp} + [A_i^* g^2(|T_i^*|)A_i]^{rp} \Big\|^{\frac{1}{r}} - \inf_{\|x\|=1} \eta(x),$$

where $T_i, A_i, B_i \in \mathbb{B}(\mathscr{H})$ ($1 \leq i \leq n$), $f$ and $g$ are nonnegative continuous functions on $[0,\infty)$ satisfying $f(t)g(t) = t$ for all $t \in [0,\infty)$, $p, r \geq 1$, $N \in \mathbb{N}$ and

$$\eta(x) = \frac{1}{2} \sum_{i=1}^n \sum_{j=1}^N \Big( \sqrt[2^j]{\langle (A_i^* g^2(|T_i^*|)A_i)^p x, x \rangle^{2^{j-1}-k_j} \langle (B_i^* f^2(|T_i|)B_i)^p x, x \rangle^{k_j}}$$
$$- \sqrt[2^j]{\langle (B_i^* f^2(|T_i|)B_i)^p x, x \rangle^{k_j+1} \langle (A_i^* g^2(|T_i^*|)A_i)^p x, x \rangle^{2^{j-1}-k_j-1}} \Big)^2.$$

## 1. INTRODUCTION

Let $\mathbb{B}(\mathscr{H})$ denote the $C^*$-algebra of all bounded linear operators on a complex Hilbert space $\mathscr{H}$ with an inner product $\langle \cdot, \cdot \rangle$ and the corresponding norm $\| \cdot \|$. In the case when $\dim \mathscr{H} = n$, we identify $\mathbb{B}(\mathscr{H})$ with the matrix algebra $\mathbb{M}_n$ of all $n \times n$ matrices with entries in the complex field. The numerical radius of $T \in \mathbb{B}(\mathscr{H})$ is defined by

$$w(T) := \sup\{|\langle Tx, x\rangle| : x \in \mathscr{H}, \| x \| = 1\}.$$

It is well known that $w(\cdot)$ defines a norm on $\mathbb{B}(\mathscr{H})$, which is equivalent to the usual operator norm $\| \cdot \|$. In fact, for any $T \in \mathbb{B}(\mathscr{H})$, $\frac{1}{2}\|T\| \leq w(T) \leq \|T\|$; see [6]. The quantity $w(T)$ is useful in studying perturbation, convergence and approximation problems as well as interactive method, etc. For more information see [1, 2, 4, 5, 7, 8, 13, 19] and references therein.

The classical Young inequality says that if $0 \leq \nu \leq 1$, then $a^\nu b^{1-\nu} \leq \nu a + (1 -$

---

$\nu)b$ $(a, b > 0)$. During the last decades several generalizations, reverses, refinements and applications of the Young inequality in various settings have been given (see [3, 12] and references therein). A refinement of the scalar Young inequality is presented in [12] as follows:

$$a^{\nu}b^{1-\nu} \leq \nu a + (1-\nu)b - r_0(a^{\frac{1}{2}} - b^{\frac{1}{2}})^2, \qquad (1.1)$$

where $r_0 = \min\{\nu, 1-\nu\}$.

Recently, Sababheh and Choi in [15] obtained a refinement of the Young inequality

$$a^{\nu}b^{1-\nu} \leq \nu a + (1-\nu)b - S_N(\nu), \qquad (1.2)$$

in which

$$S_N(\nu) := \sum_{j=1}^{N} \left( (-1)^{r_j} 2^{j-1}\nu + (-1)^{r_j+1} \left[\frac{r_j+1}{2}\right] \right) \left( \sqrt[2^j]{b^{2^{j-1}-k_j} a^{k_j}} - \sqrt[2^j]{a^{k_j+1} b^{2^{j-1}-k_j-1}} \right)^2,$$

where $N \in \mathbb{N}$, $r_j = [2^j \nu]$ and $k_j = [2^{j-1}\nu]$. Here $[x]$ is the greatest integer less than or equal to $x$. When $N = 1$, inequality (1.2) reduces to (1.1).

It follows from $\nu a + (1-\nu)b \leq (\nu a^r + (1-\nu)b^r)^{\frac{1}{r}}$ $(r \geq 1)$ and inequality (1.1) that

$$a^{\nu}b^{1-\nu} \leq (\nu a^r + (1-\nu)b^r)^{\frac{1}{r}} - S_N(\nu).$$

In particular, for $\nu = \frac{1}{2}$ we get

$$a^{\frac{1}{2}}b^{\frac{1}{2}} \leq (\frac{1}{2})^{\frac{1}{r}}(a^r + b^r)^{\frac{1}{r}} - \frac{1}{2}\sum_{j=1}^{N} \left( \sqrt[2^j]{b^{2^{j-1}-k_j} a^{k_j}} - \sqrt[2^j]{a^{k_j+1} b^{2^{j-1}-k_j-1}} \right)^2.$$

If $N = 1$, then we reach to inequality (2.1) in [12] as follows:

$$a^{\frac{1}{2}}b^{\frac{1}{2}} \leq (\frac{1}{2})^{\frac{1}{r}}(a^r + b^r)^{\frac{1}{r}} - \frac{1}{2}(a^{\frac{1}{2}} - b^{\frac{1}{2}})^2.$$

Let $T_i \in \mathbb{B}(\mathscr{H})$ $(1 \leq i \leq n)$. The Euclidean operator radius of $T_1, ..., T_n$ is defined in [14] by

$$w_e(T_1, ..., T_n) := \sup_{\|x\|=1} \left( \sum_{i=1}^{n} |\langle T_i x, x \rangle|^2 \right)^{\frac{1}{2}}.$$

In [16], the functional $w_p$ of operators $T_1, ..., T_n$ for $p \geq 1$ is defined by

$$w_p(T_1, ..., T_n) := \sup_{\|x\|=1} \left( \sum_{i=1}^{n} |\langle T_i x, x \rangle|^p \right)^{\frac{1}{p}}.$$

Let $T_1, ..., T_n \in \mathbb{B}(\mathscr{H})$. Recently, Sheikhhosseini et al. in [18] showed

$$w_p^p(A_1^* T_1 B_1, ..., A_n^* T_n B_n) \leq \frac{n^{1-\frac{1}{r}}}{2^{\frac{1}{r}}} \left\| \sum_{i=1}^{n} [B_i^* f^2(|T_i|)B_i]^{rp} + [A_i^* g^2(|T_i^*|)A_i]^{rp}) \right\|^{\frac{1}{r}} - \inf_{\|x\|=1} \zeta(x), \qquad (1.3)$$

where $\zeta(x) = \frac{1}{2}\sum_{i=1}^{n}\left(\langle [B_i^* f^2(|T_i|)B_i]^p x, x\rangle^{\frac{1}{2}} - \langle [A_i^* f^2(|T_i^*|)A_i]^p x, x\rangle^{\frac{1}{2}}\right)^2$. They also presented the following inequality

$$w_p(T_1,...,T_n) \leq \frac{1}{2}\left[\sum_{i=1}^{n}\left(\left\||T_i|^{2\alpha} + |T_i^*|^{2(1-\alpha)}\right\| - 2\inf_{\|x\|=1}\zeta_i(x)\right)^p\right]^{\frac{1}{p}}, \quad (1.4)$$

in which $0 \leq \alpha \leq 1$, $p \geq 1$ and $\zeta_i(x) = \frac{1}{2}(\langle |T_i|^{2\alpha}x, x\rangle^{\frac{1}{2}} - \langle |T_i^*|^{2(1-\alpha)}x, x\rangle^{\frac{1}{2}})^2$.

In the same paper, they showed

$$w_p^p(T_1,...,T_n) \leq \frac{1}{2}\left\|\sum_{i=1}^{n}(|T_i|^{2\alpha p} + |T_i^*|^{2(1-\alpha)p})\right\| - \inf_{\|x\|=1}\zeta(x), \quad (1.5)$$

where $\zeta(x) = \frac{1}{2}\sum_{i=1}^{n}(\langle |T_i|^{2\alpha p}x, x\rangle^{\frac{1}{2}} - \langle |T_i^*|^{2(1-\alpha)p}x, x\rangle^{\frac{1}{2}})^2$.

Moreover, they established the inequalities

$$w_p^p(T_1,...,T_n) \leq \left\|\sum_{i=1}^{n}\alpha|T_i|^p + (1-\alpha)|T_i^*|^p\right\| - \inf_{\|x\|=1}\zeta(x), \quad (1.6)$$

and

$$w_p^r(|T_1|,...,|T_n|)w_q^r(|T_1^*|,...,|T_n^*|) \leq \frac{r}{p}\|\sum_{i=1}^{n}|T_i|^p\| + \frac{r}{q}\|\sum_{i=1}^{n}|T_i^*|^q\| - \inf_{\|x\|=\|y\|=1}\delta(x,y), \quad (1.7)$$

where $\zeta(x) = \min\{\alpha, 1-\alpha\}\sum_{i=1}^{n}(\langle |T_i|^p x, x\rangle^{\frac{1}{2}} - \langle |T_i^*|^p x, x\rangle^{\frac{1}{2}})^2$ and

$$\delta(x,y) = \frac{r}{p}\left(\sqrt[2]{\sum_{i=1}^{n}\langle |T_i|x, y\rangle^p} - \sqrt[2]{\sum_{i=1}^{n}\langle |T_i^*|x, y\rangle^q}\right)^2.$$

Assume that $X \in \mathbb{B}(\mathcal{H})$. The mixed Heinz means are defined by

$$H_\alpha(A, B) = \frac{A^\alpha X B^{1-\alpha} + A^{1-\alpha} X B^\alpha}{2},$$

in which $0 \leq \alpha \leq 1$ and $A, B \geq 0$, see [10]. In [17], the authors showed that

$$w^r(A^\alpha X B^{1-\alpha}) \leq \|X\|^r \|\alpha A^r + (1-\alpha)B^r\|, \quad (1.8)$$

where $A, B, X \in \mathbb{B}(\mathcal{H})$ such that $A, B$ are positive, $r \geq 2$ and $0 \leq \alpha \leq 1$.

Using inequality (1.8), they presented an upper bound for Heinz means of matrices as follows:

$$w^r(H_\alpha(A, B)) \leq \|X\|^r \|\frac{A^r + B^r}{2}\|. \quad (1.9)$$

In this present paper, we refine inequalities (1.3)-(1.9). We also find an upper bound for the functional $w_p$.

## 2. MAIN RESULTS

To prove our numerical radius inequalities, we need several known lemmas. The first lemma is a simple result of the classical Jensen, Young and a genaralized mixed Cauchy-Schwarz inequalities [9, 11].

**Lemma 2.1.** *Let $a, b \geq 0$, $0 \leq \nu \leq 1$ and $r \neq 0$. Then*
(a) $a^\nu b^{1-\nu} \leq \nu a + (1-\nu)b \leq (\nu a^r + (1-\nu)b^r)^{\frac{1}{r}}$ *for $r \geq 1$.*
(b) *If $T \in \mathbb{B}(\mathscr{H})$ and $x, y \in \mathscr{H}$ be any vectors, then*

$$|\langle Tx, y \rangle|^2 \leq \langle |T|^{2\nu} x, x \rangle \langle |T^*|^{2(1-\nu)} y, y \rangle.$$

(c) *If $f$, $g$ are nonnegative continuous functions on $[0, \infty)$ which are satisfying the relation $f(t)g(t) = t$ $(t \in [0, \infty))$, then*

$$|\langle Tx, y \rangle| \leq \| f(|T|)x \| \| g(|T^*|)x \|$$

*for all $x, y \in \mathscr{H}$.*

**Lemma 2.2.** (*McCarty inequality* [11]). *Let $T \in \mathbb{B}(\mathscr{H})$, $T \geq 0$ and $x \in \mathscr{H}$ be a unit vector. Then*
(a) $\langle Tx, x \rangle^r \leq \langle T^r x, x \rangle$ *for $r \geq 1$;*
(b) $\langle T^r x, x \rangle \leq \langle Tx, x \rangle^r$ *for $0 < r \leq 1$.*

Now, by using inequality (1.2) we get the first result.

**Theorem 2.3.** *Let $A, B, X \in \mathbb{B}(\mathscr{H})$ such that $A, B$ are positive, $r \geq 2$ and $0 \leq \nu \leq 1$. Then*

$$w^r(A^\nu X B^{1-\nu}) \leq \|X\|^r \left[ \|\nu A^r + (1-\nu)B^r\| - \inf_{\|x\|=1} \eta(x) \right], \qquad (2.1)$$

*where*

$$\eta(x) = \sum_{j=1}^{N} \left( (-1)^{r_j} 2^{j-1} \nu + (-1)^{r_j+1} \left[ \frac{r_j+1}{2} \right] \right)$$
$$\times \left( \sqrt[2^j]{\langle B^r x, x \rangle^{2^{j-1}-k_j} \langle A^r x, x \rangle^{k_j}} - \sqrt[2^j]{\langle A^r x, x \rangle^{k_j+1} \langle B^r x, x \rangle^{2^{j-1}-k_j-1}} \right)^2.$$

*Proof.* Let $x \in \mathscr{H}$ be unit vector. Then

$$\begin{aligned}
|\langle A^\nu X B^{1-\nu} x, x\rangle|^r &= |\langle X B^{1-\nu} x, A^\nu x\rangle|^r \\
&\leq \|X\|^r \|B^{1-\nu} x\|^r \|A^\nu x\|^r \\
&= \|X\|^r \langle B^{2(1-\nu)} x, x\rangle^{\frac{r}{2}} \langle A^{2\nu} x, x\rangle^{\frac{r}{2}} \\
&\leq \|X\|^r \langle A^r x, x\rangle^\nu \langle B^r x, x\rangle^{1-\nu} \qquad \text{( by Lemma 2.2)} \\
&\leq \|X\|^r [\nu \langle A^r x, x\rangle + (1-\nu) \langle B^r x, x\rangle] \\
&\quad - \|X\|^r \sum_{j=1}^{N} \left( (-1)^{r_j} 2^{j-1} \nu + (-1)^{r_j+1} \left[\frac{r_j+1}{2}\right]\right) \times \\
&\quad \times \left( \sqrt[2^j]{\langle B^r x, x\rangle^{2^{j-1}-k_j} \langle A^r x, x\rangle^{k_j}} - \sqrt[2^j]{\langle A^r x, x\rangle^{k_j+1} \langle B^r x, x\rangle^{2^{j-1}-k_j-1}}\right)^2 \\
&\qquad\qquad\qquad\qquad\qquad\qquad\qquad \text{( by inequality (1.2))}.
\end{aligned}$$

Taking the supremum over $x \in \mathscr{H}$ with $\|x\| = 1$ in the above inequality we deduce the desired inequality. $\square$

*Remark* 2.4. Let $N = 1$ in inequality (2.1). Then

$$w^r(A^\nu X B^{1-\nu}) \leq \|X\|^r \left[ \|\nu A^r + (1-\nu) B^r\| - \inf_{\|x\|=1} \eta(x)\right], \qquad (2.2)$$

in which $\eta(x) = r_0 \left(\langle A^r x, x\rangle^{\frac{1}{2}} - \langle B^r x, x\rangle^{\frac{1}{2}}\right)^2$ and $r_0 = \min\{\nu, 1-\nu\}$. Hence inequality (2.2) is a refinement of inequality (1.8).

Using Theorem 2.3 we can find an upper bound for Heinz means of matrices that it is a refinement of (1.9).

**Theorem 2.5.** *Suppose $A, B, X \in \mathbb{B}(\mathscr{H})$ such that $A, B$ are positive. Then*

$$w^r(H_\nu(A, B)) \leq \|X\|^r \left[ \left\|\frac{A^r + B^r}{2}\right\| - \frac{1}{2} \inf \zeta(x)\right],$$

*where $r \geq 2$, $0 \leq \nu \leq 1$, $n \in \mathbb{N}$ and*

$$\zeta(x) = \sum_{j=1}^{N} \left( (-1)^{r_j} 2^{j-1} + (-1)^{r_j+1} \left[\frac{r_j+1}{2}\right]\right)$$
$$\times \left( \sqrt[2^j]{\langle B^r x, x\rangle^{2^{j-1}-k_j} \langle A^r x, x\rangle^{k_j}} - \sqrt[2^j]{\langle A^r x, x\rangle^{k_j+1} \langle B^r x, x\rangle^{2^{j-1}-k_j-1}}\right)^2.$$

*Proof.* For unit vector $x \in \mathcal{H}$, we have

$$\left|\left\langle \frac{A^\nu X B^{1-\nu} + A^{1-\nu} X B^\nu}{2} x, x \right\rangle\right|^r$$
$$\leq \left(\frac{|\langle A^\nu X B^{1-\nu} x, x\rangle| + |\langle A^{1-\nu} X B^\nu x, x\rangle|}{2}\right)^r$$
$$\leq \frac{|\langle A^\nu X B^{1-\nu} x, x\rangle|^r + |\langle A^{1-\nu} X B^\nu x, x\rangle|^r}{2}$$
$$\leq \frac{\|X\|^r}{2}[\langle \nu A^r + (1-\nu) B^r x, x\rangle] - \frac{\|X\|^r}{2} \sum_{j=1}^{N} \left((-1)^{r_j} 2^{j-1}\nu + (-1)^{r_j+1}\left[\frac{r_j+1}{2}\right]\right)$$
$$\times \left(\sqrt[2^j]{b^{2^{j-1}-k_j} a^{k_j}} - \sqrt[2^j]{a^{k_j+1} b^{2^{j-1}-k_j-1}}\right)^2$$
$$+ \frac{\|X\|^r}{2}[\langle (1-\nu)A^r + \nu B^r x, x\rangle] - \frac{\|X\|^r}{2} \sum_{j=1}^{N}\left((-1)^{r_j}2^{j-1}\nu + (-1)^{r_j+1}\left[\frac{r_j+1}{2}\right]\right)$$
$$\times \left(\sqrt[2^j]{\langle B^r x, x\rangle^{2^{j-1}-k_j}\langle A^r x, x\rangle^{k_j}} - \sqrt[2^j]{\langle A^r x, x\rangle^{k_j+1}\langle B^r x, x\rangle^{2^{j-1}-k_j-1}}\right)^2$$
$$= \|X\|^r \left[\left\langle \frac{A^r + B^r}{2} x, x\right\rangle\right] - \frac{\|X\|^r}{2} \sum_{j=1}^{N}\left((-1)^{r_j} 2^{j-1} + (-1)^{r_j+1}\left[\frac{r_j+1}{2}\right]\right)$$
$$\times \left(\sqrt[2^j]{\langle B^r x, x\rangle^{2^{j-1}-k_j}\langle A^r x, x\rangle^{k_j}} - \sqrt[2^j]{\langle A^r x, x\rangle^{k_j+1}\langle B^r x, x\rangle^{2^{j-1}-k_j-1}}\right)^2.$$

If we take the supremum over $x \in \mathcal{H}$ with $\|x\| = 1$, then we deduce the desired inequality. □

In the next theorem we show a refinement of inequality (1.3).

**Theorem 2.6.** *Let $T_i, A_i, B_i \in \mathbb{B}(\mathcal{H})$ $(1 \leq i \leq n)$ and let $f$ and $g$ be nonnegative continuous functions on $[0, \infty)$ satisfying $f(t)g(t) = t$ for all $t \in [0, \infty)$. Then*

$$w_p^p(A_1^* T_1 B_1, ..., A_n^* T_n B_n) \leq \frac{n^{1-\frac{1}{r}}}{2^{\frac{1}{r}}} \left\|\sum_{i=1}^{n}[B_i^* f^2(|T_i|)B_i]^{rp} + [A_i^* g^2(|T_i^*|)A_i]^{rp}\right\|^{\frac{1}{r}} - \inf_{\|x\|=1} \eta(x),$$
(2.3)

*where $p, r \geq 1$, $N \in \mathbb{N}$ and*

$$\eta(x) = \frac{1}{2} \sum_{i=1}^{n} \sum_{j=1}^{N} \left(\sqrt[2^j]{\langle (A_i^* g^2(|T_i^*|)A_i)^p x, x\rangle^{2^{j-1}-k_j}\langle (B_i^* f^2(|T_i|)B_i)^p x, x\rangle^{k_j}}\right.$$
$$\left. - \sqrt[2^j]{\langle (B_i^* f^2(|T_i|)B_i)^p x, x\rangle^{k_j+1}\langle (A_i^* g^2(|T_i^*|)A_i)^p x, x\rangle^{2^{j-1}-k_j-1}}\right)^2.$$

*Proof.* Let $x \in \mathscr{H}$ be any unit vector. Then

$$\sum_{i=1}^{n} |\langle A_i^* T_i B_i x, x\rangle|^p = \sum_{i=1}^{n} |\langle T_i B_i x, A_i x\rangle|^p$$

$$\leq \sum_{i=1}^{n} \|f(|T_i|)B_i x\|^p \|g(|T_i^*|)A_i x\|^p$$

(by Lemma (2.1), (c))

$$= \sum_{i=1}^{n} \langle f(|T_i|)B_i x, f(|T_i|)B_i x\rangle^{\frac{p}{2}} \langle g(|T_i^*|)A_i x, g(|T_i^*|)A_i x\rangle^{\frac{p}{2}}$$

$$= \sum_{i=1}^{n} \langle B_i^* f^2(|T_i|)B_i x, x\rangle^{\frac{p}{2}} \langle A_i^* g^2(|T_i^*|)A_i x, x\rangle^{\frac{p}{2}}$$

$$\leq \sum_{i=1}^{n} \langle (B_i^* f^2(|T_i|)B_i)^p x, x\rangle^{\frac{1}{2}} \langle (A_i^* g^2(|T_i^*|)A_i)^p x, x\rangle^{\frac{1}{2}}$$

(by Lemma (2.2), (a))

$$\leq \sum_{i=1}^{n} \left[\left(\frac{1}{2}\langle (B_i^* f^2(|T_i|)B_i)^{pr} x, x\rangle + \frac{1}{2}\langle (A_i^* g^2(|T_i^*|)A_i)^{pr} x, x\rangle\right)^{\frac{1}{r}}\right]$$

(by (1.2))

$$- \frac{1}{2}\sum_{i=1}^{n}\sum_{j=1}^{N}\Big(\sqrt[2^j]{\langle (A_i^* g^2(|T_i^*|)A_i)^p x, x\rangle^{2^{j-1}-k_j} \langle (B_i^* f^2(|T_i|)B_i)^p x, x\rangle^{k_j}}$$

$$- \sqrt[2^j]{\langle (B_i^* f^2(|T_i|)B_i)^p x, x\rangle^{k_j+1} \langle (A_i^* g^2(|T_i^*|)A_i)^p x, x\rangle^{2^{j-1}-k_j-1}}\Big)^2$$

$$\leq \frac{n^{1-\frac{1}{r}}}{2^{\frac{1}{r}}}\left\langle \left(\sum_{i=1}^{n}\left([B_i^* f^2(|T_i|)B_i]^{rp} + [A_i^* g^2(|T_i^*|)A_i]^{rp}\right)\right)x, x\right\rangle^{\frac{1}{r}}$$

$$- \frac{1}{2}\sum_{i=1}^{n}\sum_{j=1}^{N}\Big(\sqrt[2^j]{\langle (A_i^* g^2(|T_i^*|)A_i)^p x, x\rangle^{2^{j-1}-k_j} \langle (B_i^* f^2(|T_i|)B_i)^p x, x\rangle^{k_j}}$$

$$- \sqrt[2^j]{\langle (B_i^* f^2(|T_i|)B_i)^p x, x\rangle^{k_j+1} \langle (A_i^* g^2(|T_i^*|)A_i)^p x, x\rangle^{2^{j-1}-k_j-1}}\Big)^2$$

By taking supremum on unit vector $x$ in $\mathscr{H}$ we reach the desired inequality. □

**Corollary 2.7.** *Let $A_i, B_i \in \mathbb{B}(\mathscr{H})$ $(1 \leq i \leq n)$. Then for $r, p \geq 1$ we have*

$$w_p^p(A_1^* B_1, ..., A_n^* B_n) \leq \frac{n^{1-\frac{1}{r}}}{2^{\frac{1}{r}}} \left\|\sum_{i=1}^{n}\left(|B_i|^{2rp} + |A_i|^{2rp}\right)\right\|^{\frac{1}{r}} - \inf_{\|x\|=1} \eta(x),$$

*where*

$$\eta(x) = \frac{1}{2}\sum_{i=1}^{n}\sum_{j=1}^{N}\Big(\sqrt[2^j]{\langle |A_i|^{2p} x, x\rangle^{2^{j-1}-k_j} \langle |B_i|^{2p} x, x\rangle^{k_j}}$$

$$- \sqrt[2^j]{\langle |B_i|^{2p} x, x\rangle^{k_j+1} \langle |A_i|^{2p} x, x\rangle^{2^{j-1}-k_j-1}}\Big)^2.$$

*Proof.* Choosing $f(t) = g(t) = t^{\frac{1}{2}}$ and $T_i = I$ for $i = 1, 2, ..., n$ in Theorem 2.6, we get the desired result. □

**Corollary 2.8.** *Let $T_i \in \mathbb{B}(\mathscr{H})$ $(1 \leq i \leq n)$, let $f$ and $g$ be nonnegative continuous functions on $[0, \infty)$ such that $f(t)g(t) = t$ for all $t \in [0, \infty)$ and $r, p \geq 1$. Then*

$$w_p^p(T_1, ..., T_n) \leq \frac{n^{1-\frac{1}{r}}}{2^{\frac{1}{r}}} \left\| \sum_{i=1}^{n} \left( |f^{2rp}(|T_i|) + g^{2rp}(|T_i^*|) \right) \right\|^{\frac{1}{r}} - \inf_{\|x\|=1} \eta(x), \quad (2.4)$$

*where*

$$\eta(x) = \frac{1}{2} \sum_{i=1}^{n} \sum_{j=1}^{N} \Big( \sqrt[2^j]{\langle g^{2p}(|T_i^*|)x, x \rangle^{2^{j-1}-k_j} \langle f^{2p}(|T_i|)x, x \rangle^{k_j}} \\ - \sqrt[2^j]{\langle f^{2p}(|T_i|)x, x \rangle^{k_j+1} \langle g^{2p}(|T_i^*|)x, x \rangle^{2^{j-1}-k_j-1}} \Big)^2.$$

*In particular,*

$$w_p^p(T_1, ..., T_n) \leq \frac{1}{2} \left\| \sum_{i=1}^{n} (|T_i|^{2\alpha p} + |T_i^*|^{2(1-\alpha)p}) \right\| - \inf_{\|x\|=1} \eta(x), \quad (2.5)$$

*where $0 \leq \alpha \leq 1$ and*

$$\eta(x) = \frac{1}{2} \sum_{i=1}^{n} \sum_{j=1}^{N} \Big( \sqrt[2^j]{\langle |T_i^*|^{2(1-\alpha)p}x, x \rangle^{2^{j-1}-k_j} \langle |T_i|^{2\alpha p}x, x \rangle^{k_j}} \\ - \sqrt[2^j]{\langle |T_i|^{2\alpha p}x, x \rangle^{k_j+1} \langle |T_i^*|^{2(1-\alpha)p}x, x \rangle^{2^{j-1}-k_j-1}} \Big)^2.$$

*Proof.* Selecting $A_i = B_i = I$ for $i = 1, 2, ..., n$ in Theorem 2.6, we get the first result. Letting $f(t) = t^\alpha$, $g(t) = t^{1-\alpha}$, $r = 1$ and $B_i = A_i = I$ for $i = 1, 2, ..., n$ in inequality (2.4), we reach the second inequality. □

*Remark* 2.9. Note that inequality (2.5) is a refinement of inequality (1.5), since

$$\frac{1}{2} \sum_{i=1}^{n} (\langle |T_i|^{2\alpha p}x, x \rangle^{\frac{1}{2}} - \langle |T_i^*|^{2(1-\alpha)p}x, x \rangle^{\frac{1}{2}})^2 \leq \frac{1}{2} \sum_{i=1}^{n} \sum_{j=1}^{N} \Big( \sqrt[2^j]{\langle |T_i^*|^{2(1-\alpha)p}x, x \rangle^{2^{j-1}-k_j} \langle |T_i|^{2\alpha p}x, x \rangle^{k_j}} \\ - \sqrt[2^j]{\langle |T_i|^{2\alpha p}x, x \rangle^{k_j+1} \langle |T_i^*|^{2(1-\alpha)p}x, x \rangle^{2^{j-1}-k_j-1}} \Big)^2.$$

Now by letting $n = 2$, $N = 1$, $T_1 = B$ and $T_2 = C$ in Theorem 2.6, we obtain the following consequence.

**Corollary 2.10.** *Let $B, C \in \mathbb{B}(\mathscr{H})$. Then for all $p \geq 1$ and $0 \leq \alpha \leq 1$*

$$w_p^p(B, C) \leq \frac{1}{2} \left\| |B|^{2\alpha p} + |B^*|^{2(1-\alpha)p} + |C|^{2\alpha p} + |C^*|^{2(1-\alpha)p} \right\| - \inf_{\|x\|=1} \eta(x),$$

*where*

$$\eta(x) = \frac{1}{2} \Big[ (\langle |B|^{2\alpha p}x, x \rangle^{\frac{1}{2}} - \langle |B^*|^{2(1-\alpha)p}x, x \rangle^{\frac{1}{2}})^2 + (\langle |C|^{2\alpha p}x, x \rangle^{\frac{1}{2}} - \langle |C^*|^{2(1-\alpha)p}x, x \rangle^{\frac{1}{2}})^2 \Big].$$

**Theorem 2.11.** *Let $T_i \in \mathbb{B}(\mathscr{H})$ $(1 \leq i \leq n)$. Then*

$$w_p(T_1, ..., T_n) \leq \frac{1}{2}\left[\sum_{i=}^{n}\left(\left\||T_i|^{2\alpha} + |T_i^*|^{2(1-\alpha)}\right\| - 2\inf_{\|x\|=1}\eta_i(x)\right)^p\right]^{\frac{1}{p}}, \quad (2.6)$$

*where $p \geq 1$, $0 \leq \alpha \leq 1$ and*

$$\eta_i(x) = \frac{1}{2}\sum_{j=1}^{N}\Big(\sqrt[2^j]{\langle|T_i^*|^{2(1-\alpha)p}x, x\rangle^{2^{j-1}-k_j}\langle|T_i|^{2\alpha p}x, x\rangle^{k_j}} \\ - \sqrt[2^j]{\langle|T_i|^{2\alpha p}x, x\rangle^{k_j+1}\langle|T_i^*|^{2(1-\alpha)p}x, x\rangle^{2^{j-1}-k_j-1}}\Big)^2.$$

*Proof.* By using of Lemma 2.1 and inequality (1.2), for any unit vector $x \in \mathscr{H}$ we have

$$\sum_{i=1}^{n}|\langle T_i x, x\rangle|^p \leq \sum_{i=1}^{n}(\langle|T_i|^{2\alpha}x,x\rangle^{\frac{1}{2}}\langle|T_i^*|^{2(1-\alpha)}x,x\rangle^{\frac{1}{2}})^p$$

$$\text{(by Lemma (2.1), (b))}$$

$$\leq \frac{1}{2^p}\sum_{i=1}^{n}\Big[\langle|T_i|^{2\alpha}x,x\rangle + \langle|T_i^*|^{2(1-\alpha)}x,x\rangle - \\ - \sum_{j=1}^{N}\Big(\sqrt[2^j]{\langle|T_i^*|^{2(1-\alpha)p}x,x\rangle^{2^{j-1}-k_j}\langle|T_i|^{2\alpha p}x,x\rangle^{k_j}} \\ - \sqrt[2^j]{\langle|T_i|^{2\alpha p}x,x\rangle^{k_j+1}\langle|T_i^*|^{2(1-\alpha)p}x,x\rangle^{2^{j-1}-k_j-1}}\Big)^2\Big]^p$$

$$\text{(by (1.2))}$$

$$= \frac{1}{2^p}\sum_{i=1}^{n}\Big[\langle|T_i|^{2\alpha} + |T_i^*|^{2(1-\alpha)}x,x\rangle \\ - \sum_{j=1}^{N}\Big(\sqrt[2^j]{\langle|T_i^*|^{2(1-\alpha)p}x,x\rangle^{2^{j-1}-k_j}\langle|T_i|^{2\alpha p}x,x\rangle^{k_j}} - \\ - \sqrt[2^j]{\langle|T_i|^{2\alpha p}x,x\rangle^{k_j+1}\langle|T_i^*|^{2(1-\alpha)p}x,x\rangle^{2^{j-1}-k_j-1}}\Big)^2\Big]^p.$$

Thus

$$\left(\sum_{i=1}^n |\langle T_i x, x\rangle|^p\right)^{\frac{1}{p}} \leq \frac{1}{2}\Big[\sum_{i=1}^n \Big(\langle |T_i|^{2\alpha} + |T_i^*|^{2(1-\alpha)} x, x\rangle - $$
$$-\sum_{j=1}^N \Big(\sqrt[2^j]{\langle |T_i^*|^{2(1-\alpha)p} x, x\rangle^{2^{j-1}-k_j} \langle |T_i|^{2\alpha p} x, x\rangle^{k_j}}$$
$$- \sqrt[2^j]{\langle |T_i|^{2\alpha p} x, x\rangle^{k_j+1} \langle |T_i^*|^{2(1-\alpha)p} x, x\rangle^{2^{j-1}-k_j-1}}\Big)^2\Big)^p\Big]^{\frac{1}{p}}$$
$$= \frac{1}{2}\left[\sum_{i=1}^n \Big(\langle |T_i|^{2\alpha} + |T_i^*|^{2(1-\alpha)} x, x\rangle - 2\eta_i(x)\Big)^p\right]^{\frac{1}{p}}.$$

Now, by taking the supremum over all unit vector $x \in \mathscr{H}$ we get the desired result. $\square$

*Remark* 2.12. If $N = 1$ in inequality (2.6), then we reach to inequality (1.4), it follows from

$$\frac{1}{2}(\langle |T_i|^{2\alpha} x, x\rangle^{\frac{1}{2}} - \langle |T_i^*|^{2(1-\alpha)} x, x\rangle^{\frac{1}{2}})^2 \leq \frac{1}{2}\sum_{j=1}^N \Big(\sqrt[2^j]{\langle |T_i^*|^{2(1-\alpha)p} x, x\rangle^{2^{j-1}-k_j} \langle |T_i|^{2\alpha p} x, x\rangle^{k_j}}$$
$$- \sqrt[2^j]{\langle |T_i|^{2\alpha p} x, x\rangle^{k_j+1} \langle |T_i^*|^{2(1-\alpha)p} x, x\rangle^{2^{j-1}-k_j-1}}\Big)^2,$$

that inequality (2.6) is a refinement of inequality (1.4).

**Theorem 2.13.** *Let $T_i \in \mathbb{B}(\mathscr{H})$ $(1 \leq i \leq n)$. Then for $0 \leq \alpha \leq 1$ and $p \geq 2$*

$$w_p^p(T_1, ..., T_n) \leq \left\|\sum_{i=1}^n (\alpha |T_i|^p + (1-\alpha)|T_i^*|^p)\right\| - \inf_{\|x\|=1} \eta(x), \qquad (2.7)$$

*where*

$$\eta(x) = \sum_{i=1}^n \Big(\sum_{j=1}^N \Big((-1)^{r_j} 2^{j-1}\alpha + (-1)^{r_j+1}\Big[\frac{r_j+1}{2}\Big]\Big)$$
$$\times \Big(\sqrt[2^j]{\langle |T_i^*|^p x, x\rangle^{2^{j-1}-k_j} \langle |T_i|^p x, x\rangle^{k_j}} - \sqrt[2^j]{\langle |T_i|^p x, x\rangle^{k_j+1} \langle |T_i^*|^p x, x\rangle^{2^{j-1}-k_j-1}}\Big)^2$$

*Proof.* For every unit vector $x \in \mathcal{H}$ we have

$$\sum_{i=1}^{n} |\langle T_i x, x \rangle|^p$$

$$= \sum_{i=1}^{n} (|\langle T_i x, x \rangle|^2)^{\frac{p}{2}}$$

$$\leq \sum_{i=1}^{n} (\langle |T_i|^{2\alpha} x, x \rangle \langle |T_i^*|^{2(1-\alpha)} x, x \rangle)^{\frac{p}{2}} \qquad \text{(by Lemma 2.1, (b))}$$

$$\leq \sum_{i=1}^{n} (\langle |T_i|^p x, x \rangle^{\alpha} \langle |T_i^*|^p x, x \rangle^{1-\alpha}) \qquad \text{(by Lemma 2.2, (b))}$$

$$\leq \sum_{i=1}^{n} (\alpha \langle |T_i|^p x, x \rangle + (1-\alpha) \langle |T_i^*|^p x, x \rangle) -$$

$$- \sum_{i=1}^{n} \Big( \sum_{j=1}^{N} \Big( (-1)^{r_j} 2^{j-1} \alpha + (-1)^{r_j+1} \Big[ \frac{r_j+1}{2} \Big] \Big) \Big)$$

$$\times \Big( \sqrt[2^j]{\langle |T_i^*|^p x, x \rangle^{2^{j-1}-k_j} \langle |T_i|^p x, x \rangle^{k_j}} - \sqrt[2^j]{\langle |T_i|^p x, x \rangle^{k_j+1} \langle |T_i^*|^p x, x \rangle^{2^{j-1}-k_j-1}} \Big)^2$$

$$\text{(by (1.2))}$$

$$\leq \sum_{i=1}^{n} \langle (\alpha |T_i|^p + (1-\alpha) |T_i^*|^p) x, x \rangle - \inf_{\|x\|=1} \eta(x)$$

$$= \langle \sum_{i=1}^{n} (\alpha |T_i|^p + (1-\alpha) |T_i^*|^p) x, x \rangle - \inf_{\|x\|=1} \eta(x).$$

Now by taking supremum over unit vector $x \in \mathcal{H}$ we get. □

*Remark* 2.14. If we put $N = 1$ in inequality (2.7), then we get inequality (1.6). Hence inequality (2.7) is refinement of (1.6).

In [16, Remark 3.10], the author showed

$$w_p^p(B, C) \leq \frac{1}{2} \| |B|^p + |B^*|^p + |C|^p + |C^*|^p \|, \qquad (2.8)$$

in which $B, C \in \mathbb{B}(\mathcal{H})$ and $p \geq 2$. In the following result we show a refinement of (2.8).

**Corollary 2.15.** *Let $B, C \in \mathbb{B}(\mathcal{H})$. Then for $p \geq 2$,*

$$w_p^p(B, C) \leq \frac{1}{2} \| |B|^p + |B^*|^p + |C|^p + |C^*|^p \| - \inf_{\|x\|=1} \eta(x), \qquad (2.9)$$

where $\eta(x) = \frac{1}{2}\left((\langle|B|^p x, x\rangle^{\frac{1}{2}} - \langle|B^*|^p x, x\rangle^{\frac{1}{2}})^2 - (\langle|C|^p x, x\rangle^{\frac{1}{2}} - \langle|C^*|^p x, x\rangle^{\frac{1}{2}})^2\right).$
In particular, if $A \in \mathbb{B}(\mathcal{H})$, then

$$w^2(A) \leq \frac{1}{2}\|A^*A + AA^*\|.$$

*Proof.* If we take $N = 1$, $n = 2$, $T_1 = B$, $T_2 = C$, and $\alpha = \frac{1}{2}$ in Theorem 2.13, we get the first inequality.

In particular case, let $A = B + iC$ be the Cartesian decomposition of $A$. Then $A^*A + AA^* = 2(B^2 + C^2)$, and $\inf_{\|x\|=1} \eta(x) = 0$. Thus, for $p = 2$, inequality (2.9) can be written as

$$w_2^2(B, C) \leq \|B^2 + C^2\| = \frac{1}{2}\|A^*A + AA^*\|.$$

The desired inequality follows by noting that

$$w_2^2(B, C) = \sup_{\|x\|=1}\{|\langle Bx, x\rangle|^2 + |\langle Cx, x\rangle|^2\} = \sup_{\|x\|=1} |\langle Ax, x\rangle|^2 = w^2(A).$$

□

**Theorem 2.16.** *Let $T_i \in \mathbb{B}(\mathcal{H})$ $(1 \leq i \leq n)$, $r \geq 1$, and $p \geq q \geq 1$ with $\frac{1}{p} + \frac{1}{q} = \frac{1}{r}$. Then*

$$w_p^r(|T_1|, ..., |T_n|)w_q^r(|T_1^*|, ..., |T_n^*|) \leq \frac{r}{p}\|\sum_{i=1}^n |T_i|^p\| + \frac{r}{q}\|\sum_{i=1}^n |T_i^*|^q\| - \inf_{\|x\|=\|y\|=1} \lambda(x, y), \quad (2.10)$$

*where*

$$\lambda(x, y) = \sum_{j=1}^N \left((-1)^{r_j} 2^{j-1}(\frac{r}{p}) + (-1)^{r_j+1}\left[\frac{r_j+1}{2}\right]\right)\left(\sqrt[2^j]{(\sum_{i=1}^n \langle|T_i^*|x, y\rangle^q)^{2^{j-1}-k_j}(\sum_{i=1}^n \langle|T_i|x, y\rangle^p)^{k_j}}\right.$$

$$\left. - \sqrt[2^j]{(\sum_{i=1}^n \langle|T_i|x, y\rangle^p)^{k_j+1}(\sum_{i=1}^n \langle|T_i^*|x, y\rangle^q)^{2^{j-1}-k_j-1}}\right)^2.$$

*Proof.* Let $x, y \in \mathbb{B}(\mathscr{H})$ be unit vectors. Applying inequality (1.2), we get

$$\left(\left(\sum_{i=1}^{n}\langle|T_i|x,y\rangle^p\right)^{\frac{1}{p}}\left(\sum_{i=1}^{n}\langle|T_i^*|x,y\rangle^q\right)^{\frac{1}{q}}\right)^r \leq \frac{r}{p}\sum_{i=1}^{n}\langle|T_i|x,y\rangle^p + \frac{r}{q}\sum_{i=1}^{n}\langle|T_i^*|x,y\rangle^q$$

$$-\sum_{j=1}^{N}\left((-1)^{r_j}2^{j-1}(\frac{r}{p}) + (-1)^{r_j+1}\left[\frac{r_j+1}{2}\right]\right)$$

$$\times \left(\sqrt[2^j]{(\sum_{i=1}^{n}\langle|T_i^*|x,y\rangle^q)^{2^{j-1}-k_j}(\sum_{i=1}^{n}\langle|T_i|x,y\rangle^p)^{k_j}}\right.$$

$$\left. - \sqrt[2^j]{(\sum_{i=1}^{n}\langle|T_i|x,y\rangle^p)^{k_j+1}(\sum_{i=1}^{n}\langle|T_i^*|x,y\rangle^q)^{2^{j-1}-k_j-1}}\right)^2$$

(by (1.2))

$$\leq \frac{r}{p}\sum_{i=1}^{n}\langle|T_i|^p x,y\rangle + \frac{r}{q}\sum_{i=1}^{n}\langle|T_i^*|^q x,y\rangle -$$

$$-\sum_{j=1}^{N}\left((-1)^{r_j}2^{j-1}(\frac{r}{p}) + (-1)^{r_j+1}\left[\frac{r_j+1}{2}\right]\right)$$

$$\times \left(\sqrt[2^j]{(\sum_{i=1}^{n}\langle|T_i^*|x,y\rangle^q)^{2^{j-1}-k_j}(\sum_{i=1}^{n}\langle|T_i|x,y\rangle^p)^{k_j}}\right.$$

$$\left. - \sqrt[2^j]{(\sum_{i=1}^{n}\langle|T_i|x,y\rangle^p)^{k_j+1}(\sum_{i=1}^{n}\langle|T_i^*|x,y\rangle^q)^{2^{j-1}-k_j-1}}\right)^2$$

(by Lemma (2.2), (a))

$$= \frac{r}{p}\left\langle\left(\sum_{i=1}^{n}|T_i|^p\right)x,y\right\rangle + \frac{r}{q}\left\langle\left(\sum_{i=1}^{n}|T_i^*|^q\right)x,y\right\rangle$$

$$-\sum_{j=1}^{N}\left((-1)^{r_j}2^{j-1}(\frac{r}{p}) + (-1)^{r_j+1}\left[\frac{r_j+1}{2}\right]\right)$$

$$\times \left(\sqrt[2^j]{(\sum_{i=1}^{n}\langle|T_i^*|x,y\rangle^q)^{2^{j-1}-k_j}(\sum_{i=1}^{n}\langle|T_i|x,y\rangle^p)^{k_j}}\right.$$

$$\left. - \sqrt[2^j]{(\sum_{i=1}^{n}\langle|T_i|x,y\rangle^p)^{k_j+1}(\sum_{i=1}^{n}\langle|T_i^*|x,y\rangle^q)^{2^{j-1}-k_j-1}}\right)^2$$

By taking supremum on $x, y \in \mathscr{H}$ with $\|x\| = \|y\| = 1$, we get desired inequality. $\square$

*Remark* 2.17. If we put $N = 1$ in inequality (2.10), then we reach to inequality (1.7). It follows from

$$\frac{r}{p}\Big(\sqrt[2]{\sum_{i=1}^n \langle |T_i|x,y\rangle^p} - \sqrt[2]{\sum_{i=1}^n \langle |T_i^*|x,y\rangle^q}\Big)^2$$
$$\leq \sum_{j=1}^N \Big((-1)^{r_j} 2^{j-1}(\frac{r}{p}) + (-1)^{r_j+1}\Big[\frac{r_j+1}{2}\Big]\Big)\Big(\sqrt[2^j]{(\sum_{i=1}^n \langle |T_i^*|x,y\rangle^q)^{2^{j-1}-k_j}(\sum_{i=1}^n \langle |T_i|x,y\rangle^p)^{k_j}}$$
$$- \sqrt[2^j]{(\sum_{i=1}^n \langle |T_i|x,y\rangle^p)^{k_j+1}(\sum_{i=1}^n \langle |T_i^*|x,y\rangle^q)^{2^{j-1}-k_j-1}}\Big)^2,$$

that inequality (2.10) is refinement of (1.7).

**Corollary 2.18.** *Let* $T_i \in \mathbb{B}(\mathscr{H})$ $(1 \leq i \leq n)$. *Then*

$$w_e(|T_1|,...,|T_n|)w_e(|T_1^*|,...,|T_n^*|) \leq \frac{1}{2}\Big(\|\sum_{i=1}^n T_i^* T_i\| + \|\sum_{i=1}^n T_i T_i^*\|\Big) - \inf_{\|x\|=\|y\|=1} \lambda(x,y),$$

*where*

$$\lambda(x,y) = \sum_{j=1}^N \Big((-1)^{r_j} 2^{j-1}(\frac{1}{2}) + (-1)^{r_j+1}\Big[\frac{r_j+1}{2}\Big]\Big)\Big(\sqrt[2^j]{(\sum_{i=1}^n \langle |T_i^*|x,y\rangle)^{2(2^{j-1}-k_j)}(\sum_{i=1}^n \langle |T_i|x,y\rangle)^{2k_j}}$$
$$- \sqrt[2^j]{(\sum_{i=1}^n \langle |T_i|x,y\rangle)^{2(k_j+1)}(\sum_{i=1}^n \langle |T_i^*|x,y\rangle)^{2(2^{j-1}-k_j-1)}}\Big)^2.$$

*Proof.* The result obtained by letting $p = q = 2$ and $r = 1$ in inequality (2.10). □

**Corollary 2.19.** *Let* $T_1,...,T_n \in \mathbb{B}(\mathscr{H})$ *be positive operators. Then*

$$w_e(T_1,...,Tn) \leq \|\sum_{i=1}^n T_i^2\|^{\frac{1}{2}}.$$

FURTHER REFINEMENTS OF GENERALIZED NUMERICAL RADIUS INEQUALITIES   154. M. Bakherad and Kh. Shebrawi, *Upper bounds for numerical radius inequalities involving off-diagonal operator matrices*, Ann. Funct. Anal. (2018) (to appear).
5. M. Boumazgour, H. Nabwey, *A note concerning the numerical range of a basic elementary operator*, Ann. Funct. Anal. **7** (2016), no. 3, 434–441.
6. K.E. Gustafson and D.K.M. Rao, *Numerical Range, The Field of Values of Linear Operators and Matrices*, Springer, New York, 1997.
7. M. Hajmohamadi, R. Lashkaripour, M. Bakherad, *Some generalizations of numerical radius on off-diagonal part of $2 \times 2$ operator matrices.*, J. Math. Inequal. **12**(2) (2018), 447457.
8. M. Hajmohamadi, R. Lashkaripour, M. Bakherad, *Some extensions of the Young and Heinz inequalities for matrices*, Bull. Iranian Math. Soc. (2017) (to appear).
9. G.H. Hardy, J.E. Littlewood and G. Polya, *inequalities*, 2nd ed., Cambridge Univ. Press, Cambridge, (1988).
10. R. Kaur, M.S. Moslehian, M. Singh and C. Conde, *Further refinements of the Heinz inequality*, Linear Algebra Appl. **447** (2014), 26–37.
11. F. Kittaneh, *Notes on some inequalitis for Hilbert space operators*, Publ. Res. Inst. Math. Sci. **24**(2) (1988), 283–293.
12. F. Kittaneh and Y. Manasrah, *Improved Young and Heinz inequalities for matrices*, J. Math. Anal. Appl. **361** (2010), 262–269.
13. M.S. Moslehian and M. Sattari,, *Inequalities for operator space numerical radius of 2×2 block matrices*,J. Math. Phys. **57**(1) (2016), 015201, 15 pp.
14. G. Popescu, *Unitary invariants in multivariable operator theory*, Mem. Amer. Math. Soc. **200**(941) (2009).
15. M. Sababheh and D. Choi, *A complete refinement of Yongs inequality*, J. Math. Anal. Appl. **440** (2016), 379–393.
16. M. Sattari, M.S. Moslehian and K. Shebrawi, *Extension of Euclidean operator radius inequalities*, Math. Scand. **20** (2017), 129–144.
17. M. Sattari, M.S. Moslehian and T. Yamazaki, *Some genaralized numerical radius inequalities for Hilbert space operators*, Linear Algebra Appl, **470** (2014), 1–12.
18. A. Sheikhhosseini, M.S. Moslehian and K. Shebrawi, *Inequalities for generalized Euclidean operator radius via Young's inequality*, J. Math. Anal. Appl. **445**(2) (2016), 1516–1529.
19. A. Zamani, *Some lower bounds for the numerical radius of Hilbert space operators,* Adv. Oper. Theory **2** (2) (2017), 98–107.
[1,2,3] DEPARTMENT OF MATHEMATICS, FACULTY OF MATHEMATICS, UNIVERSITY OF SISTAN AND BALUCHESTAN, ZAHEDAN, I.R.IRAN.

*E-mail address*: [1]monire.hajmohamadi@yahoo.com

*E-mail address*: [2]lashkari@hamoon.usb.ac.ir

*E-mail address*: [3]mojtaba.bakherad@yahoo.com; bakherad@member.ams.org